\documentclass[11pt]{article}
\usepackage{amssymb,amsmath,amsthm}
\def\qed{\hfill $\Box$ \smallskip}
\newtheorem{theorem}{Theorem}[section]
\newtheorem{lemma}[theorem]{Lemma}

\theoremstyle{definition}

\usepackage{esint}
\usepackage{xcolor}

\usepackage{color}
\usepackage{hyperref}
\hypersetup{
	colorlinks = true,%
	citecolor = blue,
    citebordercolor = blue,
	filecolor=red,%
	linkcolor = red,
	anchorcolor = red,
	pagecolor = red,
	urlcolor= [rgb]{0.65,0.0,0.0}
	linktocpage=true,
	pdfpagelabels=true,
	bookmarksnumbered=true,
}

\begin{document}
\title{\bf{The minimizing problem involving p--Laplacian and Hardy--Littlewood--Sobolev upper critical exponent}
 \footnote{This research was supported by National Natural Science Foundation of China 11671403.}}
\date{}
\author{
Yu Su\thanks{E-mail: yizai52@qq.com (Y. Su).}
\quad
Haibo Chen\thanks{{\small Corresponding author: E-mail: math\_chb@163.com (H. Chen).}}\\
{\small School of Mathematics and Statistics, Central South University,}\\
{\small Changsha, 410083 Hunan, P.R.China.}}

\maketitle
\begin{center}
 \begin{minipage}{14cm}
\small  {\bf Abstract:}
In this paper,
we study the minimizing problem:
$$
S_{p,1,\alpha,\mu}:=
\inf_{u\in W^{1,p}(\mathbb{R}^{N})\setminus\{0\}}
\frac{
\int_{\mathbb{R}^{N}}|\nabla u|^{p}\mathrm{d}x
-
\mu
\int_{\mathbb{R}^{N}}
\frac{|u|^{p}}{|x|^{p}}
\mathrm{d}x}
{\left(
\int_{\mathbb{R}^{N}}
\int_{\mathbb{R}^{N}}
\frac{|u(x)|^{p^{*}_{\alpha}}|u(y)|^{p^{*}_{\alpha}}}{|x-y|^{\alpha}}
\mathrm{d}x
\mathrm{d}y
\right)^{\frac{p}{2\cdot p^{*}_{\alpha}}}},
$$
where
$N\geqslant3$,
$p\in(1,N)$,
$\mu\in
\left[
0,
\left(
\frac{N-p}{p}
\right)^{p}
\right)$,
$\alpha\in(0,N)$
and
$p^{*}_{\alpha}=
\frac{p}{2}\left(\frac{2N-\alpha}{N-p}\right)$
is the Hardy--Littlewood--Sobolev upper critical exponent.
Firstly,
by using refinement of Hardy-Littlewood-Sobolev inequality,
we prove that
$S_{p,1,\alpha,\mu}$
is achieved in
$\mathbb{R}^{N}$
by a radially symmetric, nonincreasing and nonnegative function.
Secondly,
we give a  estimation of extremal function.

\end{minipage}
\end{center}
 \vskip1.5mm
\section{Introduction}
In this paper,
we consider the minimizing problem:
$$
S_{p,1,\alpha,\mu}:=
\inf_{u\in W^{1,p}(\mathbb{R}^{N})\setminus\{0\}}
\frac{
\int_{\mathbb{R}^{N}}|\nabla u|^{p}\mathrm{d}x
-
\mu
\int_{\mathbb{R}^{N}}
\frac{|u|^{p}}{|x|^{p}}
\mathrm{d}x}
{\left(
\int_{\mathbb{R}^{N}}
\int_{\mathbb{R}^{N}}
\frac{|u(x)|^{p^{*}_{\alpha}}|u(y)|^{p^{*}_{\alpha}}}{|x-y|^{\alpha}}
\mathrm{d}x
\mathrm{d}y
\right)^{\frac{p}{2\cdot p^{*}_{\alpha}}}},
\eqno(\mathcal{P})
$$
where
$N\geqslant3$,
$p\in(1,N)$,
$\mu\in
\left[
0,
\left(
\frac{N-p}{p}
\right)^{p}
\right)$,
$\alpha\in(0,N)$
and
$p^{*}_{\alpha}=
\frac{p}{2}\left(\frac{2N-\alpha}{N-p}\right)$
is the Hardy--Littlewood--Sobolev upper critical exponent.

The paper was motivated by some works appeared in recent years.
For
$p=2$,
problem
$(\mathcal{P})$
is closely related to the nonlinear Choquard equation as follows:
\begin{equation}\label{1}
-\Delta u
+
V(x)u
=
\left(
|x|^{\alpha}*|u|^{q}
\right)
|u|^{q-2}u,
\mathrm{~in~}
\mathbb{R}^{N},
\end{equation}
where
$\alpha\in(0,N)$
and
$\frac{2N-\alpha}{N}\leqslant q\leqslant\frac{2N-\alpha}{N-2}$.
For
$q=2$
and
$\alpha=1$,
the equation
(\ref{1})
goes back to the description of the quantum theory of a polaron at rest by  Pekar in 1954
\cite{Pekar1954}
and the modeling of an electron trapped in its own hole in
1976
in the work of
Choquard,
as a certain approximation to Hartree--Fock theory of one--component plasma
\cite{Penrose1996}.
For
$q=\frac{2N-1}{N-2}$
and
$\alpha=1$,
by using the Green function,
it is obvious that
equation
$(\ref{1})$
can be regarded as a generalized version of
Schr\"{o}dinger--Newton system:
\begin{equation*}
\begin{aligned}
\begin{cases}
-\Delta u
+
V(x)u
=
|u|^{\frac{N+1}{N-2}}
\phi,
&
\mathrm{~in~}
\mathbb{R}^{N},\\
-\Delta \phi
=
|u|^{\frac{N+1}{N-2}},&
\mathrm{~in~}
\mathbb{R}^{N}.
\end{cases}
\end{aligned}
\end{equation*}
The existence and qualitative properties of solutions of Choquard type equations
(\ref{1})
have been widely studied in the last decades
(see \cite{Moroz2016}).
Moroz and Van Schaftingen
\cite{Moroz2015}
considered equation (\ref{1}) with lower critical exponent
$\frac{2N-\alpha}{N}$
if the potential
$1-V(x)$
should not decay to zero at infinity faster than the inverse of
$|x|^{2}$.
In
\cite{O.Alves2016},
the authors studied the equation (\ref{1}) with critical growth in the sense of Trudinger--Moser inequality and studied the existence and concentration of the ground states.
In 2016,
Gao and Yang \cite{Gao2016} firstly investigated the following critical Choquard equation:
\begin{equation}\label{2}
\begin{aligned}
-\Delta u
=
\left(
\int_{\mathbb{R}^{N}}
\frac{|u|^{2^{*}_{\alpha}}}{|x-y|^{\alpha}}
\mathrm{d}y
\right)
|u|^{2^{*}_{\alpha}-2}u
+\lambda u,
\mathrm{~in~}
\Omega,
\end{aligned}
\end{equation}
where
$\Omega$
is a bounded domain of
$\mathbb{R}^{N}$,
with lipschitz boundary,
$N\geqslant3$,
$\alpha\in(0,N)$
and
$\lambda>0$.
By using variational methods,
they established the existence, multiplicity and nonexistence of nontrivial solutions to equation (\ref{2}).
In 2017,
Mukherjee and Sreenadh \cite{Gao2016} considered the following fractional Choquard equation:
\begin{equation}\label{3}
\begin{aligned}
(-\Delta)^{s} u
=
\left(
\int_{\mathbb{R}^{N}}
\frac{|u|^{2^{*}_{\alpha,s}}}{|x-y|^{\alpha}}
\mathrm{d}y
\right)
|u|^{2^{*}_{\alpha,s}-2}u
+\lambda u,
\mathrm{~in~}
\Omega,
\end{aligned}
\end{equation}
where
$\Omega$
is a bounded domain of
$\mathbb{R}^{N}$
with $C^{1,1}$ boundary,
$s\in(0,1)$,
$N\geqslant2s$,
$\alpha\in(0,N)$
and
$\lambda>0$,
$2^{*}_{\alpha,s}=\frac{2N-\alpha}{N-2s}$
is the critical exponent in the sense of Hardy--Littlewood--Sobolev inequality.
By using variational methods,
they established the existence, multiplicity and nonexistence of nontrivial solutions to problem (\ref{3}).

For
$p\not=2$,
in 2017,
Pucci,
Xiang
and
Zhang
\cite{Pucci2017}
studied
the
Schr\"{o}dinger--Choquard--Kirchhoff equations involving the fractional p--Laplacian as follows:
\begin{equation}\label{4}
\begin{aligned}
(a+b\|u\|_{s}^{p(\theta-1)})
[(-\Delta)^{s}_{p} u+V(x)|u|^{p-2}u]
=
\lambda
f(x,u)
+
\left(
\int_{\mathbb{R}^{N}}
\frac{|u|^{p^{*}_{\alpha,s}}}{|x-y|^{\alpha}}
\mathrm{d}y
\right)
|u|^{p^{*}_{\alpha,s}-2}u
\mathrm{~in~}
\mathbb{R}^{N},
\end{aligned}
\end{equation}
where
$\|u\|_{s}=\left(
\int_{\mathbb{R}^{N}}
\int_{\mathbb{R}^{N}}
\frac{|u(x)-u(y)|^{p}}{|x-y|^{N+ps}}
\mathrm{d}x
\mathrm{d}y
+
\int_{\mathbb{R}^{N}}
V(x)|u|^{p}
\mathrm{d}x
\right)$,
$a,b\in \mathbb{R}^{+}_{0}$
with
$a+b>0$,
$\lambda>0$
is a parameter,
$s\in (0,1)$,
$N>ps$,
$\theta\in [1,\frac{N}{N-ps})$,
$\alpha\in(0,N)$,
$p^{*}_{\alpha,s}=\frac{p(2N-\alpha)}{2(N-sp)}$
is the critical exponent in the sense of Hardy--Littlewood--Sobolev inequality,
and
$f:\mathbb{R}^{N}\rightarrow\mathbb{R}$
is a Caratheodory function,
$V:\mathbb{R}^{N}\rightarrow\mathbb{R}^{+}$
is a potential function.
By using variational methods,
they established the existence of nontrivial nonnegative solution to problem (\ref{4}).

There is an open problem in
\cite{Pucci2017}.
We define the best constant:
\begin{equation}\label{5}
\begin{aligned}
S_{p,s,\alpha,\mu}:=&
\inf_{u\in W^{s,p}(\mathbb{R}^{N})\setminus\{0\}}
\frac{\int_{\mathbb{R}^{N}}
\int_{
\mathbb{R}^{N}}
\frac{|u(x)-u(y)|^{p}}{|x-y|^{N+ps}}
\mathrm{d}x
\mathrm{d}y
-
\mu
\int_{\mathbb{R}^{N}}
\frac{|u|^{p}}{|x|^{ps}}
\mathrm{d}x}
{\left(
\int_{\mathbb{R}^{N}}
\int_{\mathbb{R}^{N}}
\frac{|u(x)|^{p^{*}_{\alpha,s}}|u(y)|^{p^{*}_{\alpha,s}}}{|x-y|^{\alpha}}
\mathrm{d}x
\mathrm{d}y
\right)^{\frac{p}{2\cdot p^{*}_{\alpha,s}}}},
\end{aligned}
\end{equation}
where
$N\geqslant3$,
$p\in(1,N)$,
$s\in(0,1]$,
$\alpha\in(0,N)$
and
$\mu\in[0,\mathcal{C}_{N,s,p})$,
$\mathcal{C}_{N,s,p}$
is defined in
\cite[Theorem 1.1]{Frank2008}.
And
$p^{*}_{\alpha,s}=\frac{p(2N-\alpha)}{2(N-sp)}$
is the critical exponent in the sense of Hardy--Littlewood--Sobolev inequality.

\noindent
{\bf  Open problem:}
Is the best constant
$S_{p,s,\alpha,\mu}$
achieved?

\noindent
{\bf (Result 1)}
For
$p=2$,
$s=1$,
$\mu=0$
and
$\alpha\in(0,N)$,
Gao and Yang \cite{Gao2016}
showed that
$S_{2,1,\alpha,0}$
is achieved in
$\mathbb{R}^{N}$
by the extremal function:
\begin{equation*}
\begin{aligned}
w_{\sigma}(x)
=
\mathfrak{C}_{1}\sigma^{-\frac{N-2}{2}}
w(x),
~~
w(x)
=
\frac{b_{1}}{(b_{1}^{2}
+
|x-a_{1}|^{2})^{\frac{N-2}{2}}},
\end{aligned}
\end{equation*}
where
$\mathfrak{C}_{1}>0$
is a fixed constant,
$a_{1}\in \mathbb{R}^{N}$
and
$b_{1}\in (0,\infty)$.

\noindent
{\bf (Result 2)}
For
$p=2$,
$s\in(0,1)$,
$\mu=0$
and
$\alpha\in(0,N)$,
Mukherjee and Sreenadh \cite{Gao2016}
proved that
$S_{2,s,\alpha,0}$
is achieved in
$\mathbb{R}^{N}$
by the extremal function:
\begin{equation*}
\begin{aligned}
w_{\sigma}(x)
=
\mathfrak{C}_{2}\sigma^{-\frac{N-2s}{2}}
w(x),
~~
w(x)
=
\frac{b_{2}}{(b_{2}^{2}
+
|x-a_{2}|^{2})^{\frac{N-2s}{2}}},
\end{aligned}
\end{equation*}
where
$\mathfrak{C}_{2}>0$
is a fixed constant,
$a_{2}\in \mathbb{R}^{N}$
and
$b_{2}\in (0,\infty)$.

\noindent
{\bf (Result 3)}
For
$p=2$,
$s\in(0,1)$,
$\mu\in
\left[
0,
4^{s}\frac{\Gamma^{2}(\frac{N+2s}{4})}{\Gamma^{2}(\frac{N-2s}{4})}
\right)$
and
$\alpha\in(0,N)$,
Yang and Wu \cite{Yang2017}
showed that
$S_{2,s,\alpha,\mu}$
is achieved in
$\mathbb{R}^{N}$.

For Open problem,
we study the case of
$p\in(1,N)$,
$s=1$,
$\mu\in
\left[
0,\left(\frac{N-p}{p}\right)^{p}
\right)$
and
$\alpha\in(0,N)$.
By using
the refinement of Sobolev inequality in
\cite[Theorem 2]{Palatucci2014},
we show that
$S_{p,1,\alpha,\mu}$
is achieved in
$\mathbb{R}^{N}$
(see Theorem \ref{theorem1}).

For the case
$p\not=2$,
one expects that the minimizers of
$S_{p,s,\alpha,\mu}$
have a form similar to the function
$\omega_{\sigma}$.
However,
it is not known the explicit formula of the extremal function.
We give the estimation of
extremal function
(see Theorem \ref{theorem2} and Theorem \ref{theorem3}).

The first main result of this paper reads as follows.
\begin{theorem}\label{theorem1}
Let
$N\geqslant3$,
$p\in(1,N)$,
$\alpha\in(0,N)$
and
$\mu\in
\left[
0,\left(\frac{N-p}{p}\right)^{p}
\right)$.
Then
$S_{p,1,\alpha,\mu}$
is achieved in
$\mathbb{R}^{N}$
by a radially symmetric, nonincreasing and nonnegative function.
\end{theorem}

The second main
result of this paper reads as follows.
For
$p=2$
and
$s\in(0,1)$,
by using
Coulomb--Sobolev space
and
endpoint refined Sobolev inequality
in
\cite{Bellazzini2016},
we give
a estimation of
extremal function.

\begin{theorem}\label{theorem2}
Let
$N\geqslant3$,
$p=2$,
$\alpha\in(0,N)$,
$s\in(0,1)$
and
$\mu\in
\left[
0,
\bar{\mu}
\right)$.
Any nonnegative minimizer
$u$
of
$S_{2,s,\alpha,\mu}$
is radially symmetric and nonincreasing,
and it satisfies for
$x\not=0$
that
\begin{align*}
C_{4}
\left(
\left(
\frac{\bar{\mu}}{\bar{\mu}-\mu}
\right)
S_{2,s,\alpha,\mu}
\right)
^{\frac{(N-\alpha)(N-2s)}{2N(N+2s-\alpha)}}
\left(
\frac{N}{\omega_{N-1}}
\right)
^{\frac{N-2s}{2N}}
\frac{1}{
|x|^{\frac{N-2s}{2}}}
\geqslant
u(x),
\end{align*}
where
$\omega_{N-1}$
is the area of the unit sphere in
$\mathbb{R}^{N}$.
\end{theorem}

The third main
result of this paper reads as follows.
For
$p\not=2$
and
$s=1$,
we give
a estimation of
extremal function.
\begin{theorem}\label{theorem3}
Let
$N\geqslant3$,
$p\in(1,N)$,
$\alpha\in(0,N)$
and
$\mu\in
\left[
0,
\tilde{\mu}
\right)$.
Any nonnegative minimizer
$u$
of
$S_{p,1,\alpha,\mu}$
is radially symmetric and nonincreasing,
and it satisfies for
$x\not=0$
that
\begin{align*}
\left(
\frac{2^{\alpha}N^{2}}{\omega^{2}_{N-1}}
\right)^{\frac{1}{2\cdot p^{*}_{\alpha}}}
\frac{1}{|x|^{\frac{N-p}{p}}}
\geqslant
u(x),
\end{align*}
where
$\omega_{N-1}$
is the area of the unit sphere in
$\mathbb{R}^{N}$.
\end{theorem}
\section{Preliminaries}
The Sobolev space
$W^{1,p}(\mathbb{R}^{N})$
is the completion of
$C^{\infty}_{0}(\mathbb{R}^{N})$
with respect to the norm
$$
\|u\|_{W}^{p}=\int_{\mathbb{R}^{N}}|\nabla u|^{p}\mathrm{d}x.
$$
For
$s\in(0,1)$
and
$p\in(1,N)$,
the fractional Sobolev space
$W^{s,p}(\mathbb{R}^{N})$
is defined by
$$
W^{s,p}(\mathbb{R}^{N}):=
\left\{
u\in L^{\frac{Np}{N-sp}}(\mathbb{R}^{N})|
\int_{\mathbb{R}^{N}}
\int_{
\mathbb{R}^{N}}
\frac{|u(x)-u(y)|^{p}}{|x-y|^{N+ps}}
\mathrm{d}x
\mathrm{d}y
<
\infty
\right\}.
$$
For
$s\in(0,1)$
and
$p\in(1,N)$,
we introduce the Hardy inequalities:
$$
\bar{\mu}
\int_{\mathbb{R}^{N}}
\frac{|u|^{2}}{|x|^{2s}}
\mathrm{d}x
\leqslant
\int_{\mathbb{R}^{N}}
\int_{
\mathbb{R}^{N}}
\frac{|u(x)-u(y)|^{2}}{|x-y|^{N+2s}}
\mathrm{d}x
\mathrm{d}y,
~\mathrm{for~any}~u\in W^{s,2}(\mathbb{R}^{N})
~\mathrm{and}~
\bar{\mu}=4^{s}\frac{\Gamma^{2}(\frac{N+2s}{4})}{\Gamma^{2}(\frac{N-2s}{4})},
$$
and
$$
\tilde{\mu}
\int_{\mathbb{R}^{N}}
\frac{|u|^{p}}{|x|^{p}}
\mathrm{d}x
\leqslant
\int_{
\mathbb{R}^{N}}
|\nabla u|^{p}
\mathrm{d}x,
~\mathrm{for~any}~u\in W^{1,p}(\mathbb{R}^{N})
~\mathrm{and}~\tilde{\mu}=\left(\frac{N-p}{p}\right)^{p}.$$
The Coulomb--Sobolev space \cite{Bellazzini2016} is defined by
\begin{align}\label{8}
\mathcal{E}^{s,\alpha,2^{*}_{\alpha,s}}(\mathbb{R}^{N})
\!
=
\!
\left\{
\!
\int_{\mathbb{R}^{N}}
\!\!
\int_{\mathbb{R}^{N}}
\!\!\!
\frac{|u(x)-u(y)|^{2}}
{|x-y|^{N+2s}}
\mathrm{d}x
\mathrm{d}y
\!<
\!
\infty
~\mathrm{and}
\int_{\mathbb{R}^{N}}
\!\!
\int_{\mathbb{R}^{N}}
\!\!\!
\frac{|u(x)|^{2^{*}_{\alpha,s}}|u(y)|^{2^{*}_{\alpha,s}}}{|x-y|^{\alpha}}
\mathrm{d}x
\mathrm{d}y<\infty
\!
\right\}.
\end{align}
We endow the space
$\mathcal{E}^{s,\alpha,2^{*}_{\alpha,s}}(\mathbb{R}^{N})$
with the norm
\begin{align}\label{9}
\|u\|_{\mathcal{E},\alpha}^{2}
=
\int_{\mathbb{R}^{N}}
\int_{\mathbb{R}^{N}}
\frac{|u(x)-u(y)|^{2}}
{|x-y|^{N+2s}}
\mathrm{d}x
\mathrm{d}y
+
\left(
\int_{\mathbb{R}^{N}}
\int_{\mathbb{R}^{N}}
\frac{|u(x)|^{2^{*}_{\alpha,s}}|u(y)|^{2^{*}_{\alpha,s}}}{|x-y|^{\alpha}}
\mathrm{d}x
\mathrm{d}y
\right)^{\frac{1}{2^{*}_{\alpha,s}}}.
\end{align}
We could define the best constant:
\begin{equation}\label{10}
\begin{aligned}
S_{p,1,0,\mu}:=&
\inf_{u\in W^{1,p}(\mathbb{R}^{N})\setminus\{0\}}
\frac{
\|u\|^{p}_{W}
-
\mu
\int_{\mathbb{R}^{N}}
\frac{|u|^{p}}{|x|^{p}}
\mathrm{d}x
}
{(\int_{\mathbb{R}^{N}}
|u|^{p^{*}}
\mathrm{d}x)^{\frac{p}{p^{*}}}},
\end{aligned}
\end{equation}
where
$S_{p,1,0,\mu}$
is attained in
$\mathbb{R}^{N}$.

\begin{lemma}\label{lemma1}
$(Hardy-Littlewood-Sobolev~inequality,~\cite{Lieb2001})$
Let
$t,r>1$
and
$0<\mu<N$
with
$\frac{1}{t}+\frac{1}{r}+\frac{\mu}{N}=2$,
$f\in L^{t}(\mathbb{R}^{N})$
and
$h\in L^{r}(\mathbb{R}^{N})$.
There exists a sharp constant
$C_{2}>0$,
independent of
$f,g$
such that
$$\int_{\mathbb{R}^{N}}\int_{\mathbb{R}^{N}}
\frac{|f(x)||h(y)|}
{|x-y|^{\mu}}
\mathrm{d}x\mathrm{d}y
\leqslant
C_{2}
\|f\|_{t}
\|h\|_{r}.$$
\end{lemma}

A measurable function
$u:\mathbb{R}^{N}\rightarrow \mathbb{R}$
belongs to the Morrey space
$\|u\|_{\mathcal{L}^{r,\varpi}(\mathbb{R}^{N})}$
with
$r\in[1,\infty)$
and
$\varpi\in(0,N]$
if and only if
$$
\|u\|
_{\mathcal{L}^{r,\varpi}(\mathbb{R}^{N})}^{r}
=
\sup_{R>0,x\in\mathbb{R}^{N}}
R^{\varpi-N}
\int_{B(x,R)}
|u(y)|^{r}
\mathrm{d}y
<\infty.
$$
\begin{lemma}\label{lemma2}
~\cite{Palatucci2014}
For any
$1< p<N$,
let
$p^{*}=\frac{Np}{N-p}$.
There exists
$C_{3}>0$
such that
for
$\theta$
and
$\vartheta$
satisfying
$\frac{p}{p^{*}}\leqslant\theta<1$,
$1\leqslant \vartheta<p^{*}=\frac{Np}{N-p}$,
we have
\begin{align*}
\left(
\int_{\mathbb{R}^{N}}
|u|^{p^{*}}
\mathrm{d}x
\right)^{\frac{1}{p^{*}}}
\leqslant
C_{3}
\|u\|_{W}^{\theta}
\|u\|_{\mathcal{L}^{\vartheta,\frac{\vartheta(N-p)}{p}}(\mathbb{R}^{N})}^{1-\theta},
\end{align*}
for any
$u\in W^{1,p}(\mathbb{R}^{N})$.
\end{lemma}
\begin{lemma}\label{lemma3}
$(Endpoint~refined~Sobolev~inequality,~
\cite[Theorem 1.2]{Bellazzini2016})$
Let
$\alpha\in(0,N)$
and
$s\in(0,1)$.
Then there exists a  constant
$C_{4}>0$
such that the inequality
\begin{align*}
\|u\|_{L^{\frac{2N}{N-2s}}(\mathbb{R}^{N})}
\leqslant&
C_{4}
\left(
\int_{R^{N}}
\int_{R^{N}}
\frac{|u(x)-u(y)|^{2}}
{|x-y|^{N+2s}}
\mathrm{d}x
\mathrm{d}y
\right)
^{\frac{(N-\alpha)(N-2s)}{2N(N+2s-\alpha)}}\\
&~~~~~~~~~\cdot
\left(
\int_{\mathbb{R}^{N}}
\int_{\mathbb{R}^{N}}
\frac{|u(x)|^{\frac{2N-\alpha}{N-2s}}|u(y)|^{\frac{2N-\alpha}{N-2s}}}{|x-y|^{\alpha}}
\mathrm{d}x
\mathrm{d}y
\right)
^{\frac{s(N-2s)}{N(N+2s-\alpha)}}
\end{align*}
holds for all
$u\in \mathcal{E}^{s,\alpha,2^{*}_{\alpha,s}}(\mathbb{R}^{N})$.
\end{lemma}
\section{The proof of Theorem \ref{theorem1}}
We show the refinement of Hardy-Littlewood-Sobolev inequality.
This inequality plays a key role in the proof of Theorem \ref{theorem1}.
\begin{lemma}\label{lemma4}
For any
$1< p<N$
and
$\alpha\in(0,N)$,
there exists
$C_{5}>0$
such that
for
$\theta$
and
$\vartheta$
satisfying
$\frac{p}{p^{*}}\leqslant\theta<1$,
$1\leqslant \vartheta<p^{*}=\frac{Np}{N-p}$,
we have
\begin{align*}
\left(
\int_{\mathbb{R}^{N}}
\int_{\mathbb{R}^{N}}
\frac{|u(x)|^{p^{*}_{\alpha}}|u(y)|^{p^{*}_{\alpha}}}{|x-y|^{\alpha}}
\mathrm{d}x
\mathrm{d}y
\right)^{\frac{1}{p^{*}_{\alpha}}}
\leqslant
C_{5}
\|u\|_{W}^{2\theta}
\|u\|_{\mathcal{L}^{\vartheta,\frac{\vartheta(N-p)}{p}}(\mathbb{R}^{N})}^{2(1-\theta)},
\end{align*}
for any
$u\in W^{1,p}(\mathbb{R}^{N})$.
\end{lemma}
\noindent
{\bf Proof.}
By using
Lemma \ref{lemma2},
we have
\begin{align}\label{11}
\left(
\int_{\mathbb{R}^{N}}
|u|^{p^{*}}
\mathrm{d}x
\right)^{\frac{1}{p^{*}}}
\leqslant
C_{3}
\|u\|_{W}^{\theta}
\|u\|_{\mathcal{L}^{\vartheta,\frac{\vartheta(N-p)}{p}}(\mathbb{R}^{N})}^{1-\theta}.
\end{align}
By Hardy-Littlewood-Sobolev inequality
and
(\ref{11}),
we obtain
\begin{equation*}
\begin{aligned}
\left(
\int_{\mathbb{R}^{N}}
\int_{\mathbb{R}^{N}}
\frac{|u(x)|^{p^{*}_{\alpha}}|u(y)|^{p^{*}_{\alpha}}}{|x-y|^{\alpha}}
\mathrm{d}x
\mathrm{d}y
\right)^{\frac{1}{p^{*}_{\alpha}}}
\leqslant&
C_{2}^{\frac{1}{p^{*}_{\alpha}}}
\|u\|_{L^{p^{*}}(\mathbb{R}^{N})}^{2}
\leqslant
C_{2}^{\frac{1}{p^{*}_{\alpha}}}
C_{3}^{2}
\|u\|_{W}^{2\theta}
\|u\|_{\mathcal{L}^{\vartheta,\frac{\vartheta(N-p)}{p}}(\mathbb{R}^{N})}^{2(1-\theta)}.
\end{aligned}
\end{equation*}
\qed

In \cite{Palatucci2014},
there is a misprint,
the authors point out it by themselves.
The right one is
\begin{align}\label{12}
L^{p^{*}}(\mathbb{R}^{N})
\hookrightarrow
\mathcal{L}^{r,r\frac{N-p}{p}}(\mathbb{R}^{N}),
\end{align}
for any
$p\in(1,N)$
and
$r\in[1,p^{*})$.
This embedding plays a key role in the proof of Theorem \ref{theorem1}.

\noindent
{\bf Proof of Theorem \ref{theorem1}:}
\noindent
{\bf Step 1.}
Suppose now
$0\leqslant\mu<\tilde{\mu}=\left(\frac{N-p}{p}\right)^{p}$.
Applying
Lemma
\ref{lemma4}
with
$\vartheta=p$,
we have
\begin{align}\label{13}
\left(
\int_{\mathbb{R}^{N}}
\int_{\mathbb{R}^{N}}
\frac{|u(x)|^{p^{*}_{\alpha}}|u(y)|^{p^{*}_{\alpha}}}{|x-y|^{\alpha}}
\mathrm{d}x
\mathrm{d}y
\right)^{\frac{1}{p^{*}_{\alpha}}}
\leqslant
C
\left(\|u\|_{W}^{p}-
\mu
\int_{\mathbb{R}^{N}}
\frac{|u|^{p}}{|x|^{p}}
\mathrm{d}x\right)^{\frac{2\theta}{p}}
\|u\|_{\mathcal{L}^{p,N-p}(\mathbb{R}^{N})}^{2(1-\theta)},
\end{align}
for
$u\in W^{1,p}(\mathbb{R}^{N})$.
Let
$\{u_{n}\}$
be a minimizing sequence of
$S_{p,1,\alpha,\mu}$,
that is
\begin{align*}
\|u_{n}\|_{W}^{p}-
\mu
\int_{\mathbb{R}^{N}}
\frac{|u_{n}|^{p}}{|x|^{p}}
\mathrm{d}x
\rightarrow
S_{p,1,\alpha,\mu},
~\mathrm{as}~
n\rightarrow\infty,
\end{align*}
and
\begin{align*}
\int_{\mathbb{R}^{N}}
\int_{\mathbb{R}^{N}}
\frac{|u_{n}(x)|^{p^{*}_{\alpha}}|u_{n}(y)|^{p^{*}_{\alpha}}}{|x-y|^{\alpha}}
\mathrm{d}x
\mathrm{d}y
=1.
\end{align*}
Inequality
(\ref{13})
enables us to find
$C>0$
independent of $n$
such that
\begin{align}\label{14}
\|u_{n}\|_{\mathcal{L}^{p,N-p}(\mathbb{R}^{N})}
\geqslant
C>0.
\end{align}
We have the chain of inclusions
\begin{align}\label{15}
W^{1,p}(\mathbb{R}^{N})
\hookrightarrow
L^{p^{*}}(\mathbb{R}^{N})
\hookrightarrow
\mathcal{L}^{p,N-p}
(\mathbb{R}^{N}),
\end{align}
which implies that
\begin{align}\label{16}
\|u_{n}\|_{\mathcal{L}^{p,N-p}(\mathbb{R}^{N})}
\leqslant
C.
\end{align}
Applying
(\ref{14})
and
(\ref{16}),
there exists
$C>0$
such that
\begin{align*}
0<C
\leqslant
\|u_{n}\|_{\mathcal{L}^{p,N-p}(\mathbb{R}^{N})}
\leqslant
C^{-1}.
\end{align*}
Combining the definition of Morrey space and above inequalities,
for all $n\in \mathbb{N}$,
we get the existence of
$\lambda_{n}>0$
and
$x_{n}\in \mathbb{R}^{N}$
such that
\begin{align*}
\frac{1}{\lambda_{n}^{p}}
\int_{B(x_{n},\lambda_{n})}
|u_{n}(y)|^{p}
\mathrm{d}y
\geqslant
\|u_{n}\|_{\mathcal{L}^{p,N-p}(\mathbb{R}^{N})}^{p}
-
\frac{C}{2n}
\geqslant
\tilde{C}
>0
\end{align*}
for some new positive constant $\tilde{C}$ that does not depend on $n$.

Let
$v_{n}(x)=\lambda_{n}^{\frac{N-p}{p}}u_{n}(\lambda_{n}x)$.
Notice that,
by using the scaling invariance,
we have
\begin{align*}
\|v_{n}\|_{W}^{p}-
\mu
\int_{\mathbb{R}^{N}}
\frac{|v_{n}|^{p}}{|x|^{p}}
\mathrm{d}x
\rightarrow
S_{p,1,\alpha,\mu},
~\mathrm{as}~
n\rightarrow\infty,
\end{align*}
and
\begin{align*}
\int_{\mathbb{R}^{N}}
\int_{\mathbb{R}^{N}}
\frac{|v_{n}(x)|^{p^{*}_{\alpha}}|v_{n}(y)|^{p^{*}_{\alpha}}}{|x-y|^{\alpha}}
\mathrm{d}x
\mathrm{d}y
=1.
\end{align*}
Then
\begin{align*}
\int_{B(\frac{x_{n}}{\lambda_{n}},1)}
|v_{n}(y)|^{p}
\mathrm{d}y
\geqslant
\tilde{C}
>0.
\end{align*}
We can
also show that
$v_{n}$
is bounded in
$W^{1,p}(\mathbb{R}^{N})$.
Hence, we may assume

$$
v_{n}\rightharpoonup v~
\mathrm{in}
~
W^{1,p}(\mathbb{R}^{N}),
v_{n}\rightarrow v
~
\mathrm{a.e. ~in}
~
\mathbb{R}^{N},
v_{n}\rightarrow v
~\mathrm{in}~
L^{q}_{loc}(\mathbb{R}^{N})~~\mathrm{for~all~}q\in[p,p^{*} ).
$$
We claim that
$\{\frac{x_{n}}{\lambda_{n}}\}$
is uniformly bounded in
$n$.
Indeed, for any $0<\beta<p$,
by H\"{o}lder's inequality,
 we observe that
\begin{align*}
0<
\tilde{C}
\leqslant&
\int_{B(\frac{x_{n}}{\lambda_{n}},1)}
|v_{n}|^{p}
\mathrm{d}y
=
\int_{B(\frac{x_{n}}{\lambda_{n}},1)}
|y|^{\frac{p\beta}{\frac{p(N-\beta)}{N-p}}}
\frac{|v_{n}|^{p}}
{|y|^{\frac{p\beta}{\frac{p(N-\beta)}{N-p}}}}
\mathrm{d}y\\
\leqslant&
\left(
\int_{B(\frac{x_{n}}{\lambda_{n}},1)}
|y|^{{\frac{\beta(N-p)}{p-\beta}}}
\mathrm{d}y
\right)
^{1-\frac{N-p}{N-\beta}}
\left(
\int_{B(\frac{x_{n}}{\lambda_{n}},1)}
\frac{|v_{n}|^{\frac{p(N-\beta)}{N-p}}}
{|y|^{\beta}}
\mathrm{d}y
\right)^{\frac{N-p}{N-\beta}}\textcolor{red}{.}
\end{align*}
By the rearrangement inequality,
see
\cite[Theorem 3.4]{Lieb2001},
we have
\begin{align*}
\int_{B(\frac{x_{n}}{\lambda_{n}},1)}
|y|^{{\frac{\beta(N-p)}{p-\beta}}}
\mathrm{d}y
\leqslant
\int_{B(0,1)}
|y|^{{\frac{\beta(N-p)}{p-\beta}}}
\mathrm{d}y
\leqslant
C.
\end{align*}
Therefore,
\begin{align}\label{17}
0<
C
\leqslant
\int_{B(\frac{x_{n}}{\lambda_{n}},1)}
\frac{|v_{n}|^{\frac{p(N-\beta)}{N-p}}}
{|y|^{\beta}}
\mathrm{d}y.
\end{align}
Now, suppose on the contrary, that
$\frac{x_{n}}{\lambda_{n}}\rightarrow\infty$
as
$n\rightarrow\infty$.
Then, for any
$y\in B(\frac{x_{n}}{\lambda_{n}},1)$,
we have
$|y|\geqslant |\frac{x_{n}}{\lambda_{n}}|-1$
for
$n$
large.
Thus,
\begin{align*}
\int_{B(\frac{x_{n}}{\lambda_{n}},1)}
\frac{|v_{n}|^{\frac{p(N-\beta)}{N-p}}}
{|y|^{\beta}}
\mathrm{d}y
\leqslant&
\frac{1}
{(|\frac{x_{n}}{\lambda_{n}}|-1)^{\beta}}
\int_{B(\frac{x_{n}}{\lambda_{n}},1)}
|v_{n}|^{\frac{p(N-\beta)}{N-p}}
\mathrm{d}y\\
\leqslant&
\frac{
\left|
B(\frac{x_{n}}{\lambda_{n}},1)
\right|
^{\frac{\beta}{N}}
}
{(|\frac{x_{n}}{\lambda_{n}}|-1)^{\beta}}
\left(
\int_{B(\frac{x_{n}}{\lambda_{n}},1)}
|v_{n}|^{\frac{Np}{N-p}}
\mathrm{d}y
\right)
^{\frac{N-\beta}{N}}\\
\leqslant&
\frac{
\left|
B(\frac{x_{n}}{\lambda_{n}},1)
\right|
^{\frac{\beta}{N}}
}
{(|\frac{x_{n}}{\lambda_{n}}|-1)^{\beta}}
\cdot
\frac{\|v_{n}\|^{\frac{N-\beta}{N}}_{W}}{S_{p,1,0,0}^{\frac{N-\beta}{N-p}}}
\leqslant
\frac{C}
{(|\frac{x_{n}}{\lambda_{n}}|-1)^{\beta}}
\rightarrow0
~\mathrm{as}~
n\rightarrow\infty,
\end{align*}
which contradicts
(\ref{17}).
Hence,
$\{\frac{x_{n}}{\lambda_{n}}\}$
is bounded,
and there exists
$R>0$
such that
\begin{align*}
\int_{B(0,R)}
|v_{n}(y)|^{p}
\mathrm{d}y
\geqslant
\int_{B(\frac{x_{n}}{\lambda_{n}},1)}
|v_{n}(y)|^{p}
\mathrm{d}y
\geqslant
\tilde{C}
>0.
\end{align*}
Since the embedding
$W^{1,p}(\mathbb{R}^{N})\hookrightarrow L_{loc}^{q}(\mathbb{R}^{N})$
$q\in[p,p^{*})$
is compact,
we deduce that
\begin{align*}
\int_{B(0,R)}
|v(y)|^{p}
\mathrm{d}y
\geqslant
\tilde{C}
>0,
\end{align*}
which means
$v\not\equiv0$.

\noindent
{\bf Step 2.}
Set
$$h(t)=t^{\frac{2\cdot p^{*}_{\alpha}}{p}},~t\geqslant0~(1<p<N).$$
Since
$p\in(1,N)$
and
$\alpha\in(0,N)$,
we get
$$\frac{2\cdot p^{*}_{\alpha}}{p}=\frac{2N-\alpha}{N-p}>1~\mathrm{and}~N+p-\alpha>0.$$
We know that
$$h^{''}(t)=\frac{(2N-\alpha)(N+p-\alpha)}{(N-p)^{2}} t^{\frac{2p-\alpha}{N-p}}\geqslant0,$$
which implies that
$h(t)$
is a convex function.
By using
$h(0)=0$
and
$l\in[0,1]$,
we know
\begin{align}\label{110}
h(lt)=h(lt+(1-t)\cdot0)\leqslant lh(t)+(1-l)h(0)=lh(t).
\end{align}
For any
$t_{1},t_{2}\in[0,\infty)$,
applying last inequality with
$l=\frac{t_{1}}{t_{1}+t_{2}}$
and
$l=\frac{t_{2}}{t_{1}+t_{2}}$,
we get
\begin{equation}\label{111}
\begin{aligned}
h(t_{1})+h(t_{2})
=&
h
\left(
(t_{1}+t_{2})\frac{t_{1}}{t_{1}+t_{2}}
\right)
+
h
\left(
(t_{1}+t_{2})\frac{t_{2}}{t_{1}+t_{2}}
\right)\\
\leqslant&
\frac{t_{1}}{t_{1}+t_{2}}
h
\left(
t_{1}+t_{2}
\right)
+
\frac{t_{2}}{t_{1}+t_{2}}
h
\left(
t_{1}+t_{2}
\right)
~(\mathrm{by}~(\ref{110}))\\
=&h
\left(
t_{1}+t_{2}
\right).
\end{aligned}
\end{equation}
Now,
we claim that
$v_{n}\rightarrow v~\mathrm{strongly~in}~W^{1,p}(\mathbb{R}^{N}).$
Set
$$K(u,v)=
\int_{\mathbb{R}^{N}}
|\nabla u|^{p-2}
\nabla u \nabla v
\mathrm{d}x
-
\mu
\int_{\mathbb{R}^{N}}
\frac{|u|^{p-2}uv}{|x|^{p}}
\mathrm{d}x.$$
Since
$\{v_{n}\}$
is a minimizing sequence,
$$\lim\limits_{n\rightarrow\infty}K(v_{n},v_{n})=S_{p,1,\alpha,\mu}.$$
By using Br\'{e}zis--Lieb type lemma
\cite{Brezis1983}
and
\cite[Theorem 2.3]{Pucci2017},
we know
\begin{align}\label{113}
K(v,v)+\lim\limits_{n\rightarrow\infty}K(v_{n}-v,v_{n}-v)=
\lim\limits_{n\rightarrow\infty}K(v_{n},v_{n})+o(1)=
S_{p,1,\alpha,\mu}+o(1),
\end{align}
and
\begin{equation}\label{114}
\begin{aligned}
&\int_{\mathbb{R}^{N}}
\int_{\mathbb{R}^{N}}
\frac{|v_{n}(x)|^{p^{*}_{\alpha}}|v_{n}(y)|^{p^{*}_{\alpha}}}{|x-y|^{\alpha}}
\mathrm{d}x
\mathrm{d}y
-
\int_{\mathbb{R}^{N}}
\int_{\mathbb{R}^{N}}
\frac{|v_{n}(x)-v(x)|^{p^{*}_{\alpha}}|v_{n}(y)-v(y)|^{p^{*}_{\alpha}}}{|x-y|^{\alpha}}
\mathrm{d}x
\mathrm{d}y\\
=&
\int_{\mathbb{R}^{N}}
\int_{\mathbb{R}^{N}}
\frac{|v(x)|^{p^{*}_{\alpha}}|v(y)|^{p^{*}_{\alpha}}}{|x-y|^{\alpha}}
\mathrm{d}x
\mathrm{d}y
+
o(1).
\end{aligned}
\end{equation}
where
$o(1)$
denotes a quantity that tends to zero as
$n\rightarrow\infty$.
Therefore,
\begin{align*}
1
=&
\lim\limits_{n\rightarrow\infty}
\int_{\mathbb{R}^{N}}
\int_{\mathbb{R}^{N}}
\frac{|v_{n}(x)|^{p^{*}_{\alpha}}|v_{n}(y)|^{p^{*}_{\alpha}}}{|x-y|^{\alpha}}
\mathrm{d}x
\mathrm{d}y\\
=&
\lim\limits_{n\rightarrow\infty}
\int_{\mathbb{R}^{N}}
\int_{\mathbb{R}^{N}}
\frac{|v_{n}(x)-v(x)|^{p^{*}_{\alpha}}|v_{n}(y)-v(y)|^{p^{*}_{\alpha}}}{|x-y|^{\alpha}}
\mathrm{d}x
\mathrm{d}y\\
&+
\int_{\mathbb{R}^{N}}
\int_{\mathbb{R}^{N}}
\frac{|v(x)|^{p^{*}_{\alpha}}|v(y)|^{p^{*}_{\alpha}}}{|x-y|^{\alpha}}
\mathrm{d}x
\mathrm{d}y\\
\leqslant&
S_{p,1,\alpha,\mu}^{-\frac{2\cdot p^{*}_{\alpha}}{p}}
\left(
\lim\limits_{n\rightarrow\infty} K(v_{n}-v,v_{n}-v)
\right)^{\frac{2\cdot p^{*}_{\alpha}}{p}}
+
S_{p,1,\alpha,\mu}^{-\frac{2\cdot p^{*}_{\alpha}}{p}}
\left(
 K(v,v)
\right)^{\frac{2\cdot p^{*}_{\alpha}}{p}}\\
\leqslant&
S_{p,1,\alpha,\mu}^{-\frac{2\cdot p^{*}_{\alpha}}{p}}
\left(
\lim\limits_{n\rightarrow\infty} K(v_{n}-v,v_{n}-v)
+
 K(v,v)
\right)^{\frac{2\cdot p^{*}_{\alpha}}{p}}
~(\mathrm{by}~(\ref{111}))\\
\leqslant&1~(\mathrm{by}~(\ref{113})).
\end{align*}
Therefore, all the inequalities above have to be equalities.
We know that
\begin{equation}\label{115}
\begin{aligned}
&
\left(
\lim\limits_{n\rightarrow\infty} K(v_{n}-v,v_{n}-v)
\right)^{\frac{2\cdot p^{*}_{\alpha}}{p}}
+
\left(
 K(v,v)
\right)^{\frac{2\cdot p^{*}_{\alpha}}{p}}\\
=&
\left(
\lim\limits_{n\rightarrow\infty} K(v_{n}-v,v_{n}-v)
+
 K(v,v)
\right)^{\frac{2\cdot p^{*}_{\alpha}}{p}}.
\end{aligned}
\end{equation}
We show that
$\lim\limits_{n\rightarrow\infty} K(v_{n}-v,v_{n}-v)=0$.
Combining
\eqref{111}
and
\eqref{115},
we know that
$$\mathrm{eithor}~\lim\limits_{n\rightarrow\infty} K(v_{n}-v,v_{n}-v)=0~\mathrm{or}~K(v,v)=0.$$
Since
$v\not\equiv0$,
so $K(v,v)\not=0$.
Therefore,
\begin{align}\label{116}
\lim\limits_{n\rightarrow\infty} K(v_{n}-v,v_{n}-v)=0.
\end{align}
This implies that
$v_{n}\rightarrow v~\mathrm{strongly~in}~W^{1,p}(\mathbb{R}^{N}).$
Moreover,
we get
\begin{align}\label{117}
\lim\limits_{n\rightarrow\infty}
\int_{\mathbb{R}^{N}}
\int_{\mathbb{R}^{N}}
\frac{|v_{n}(x)-v(x)|^{p^{*}_{\alpha}}|v_{n}(y)-v(y)|^{p^{*}_{\alpha}}}{|x-y|^{\alpha}}
\mathrm{d}x
\mathrm{d}y=0.
\end{align}

\noindent
{\bf Step 3.}
Since
$v\not\equiv0$,
putting
(\ref{116})
into
(\ref{113}),
and
inserting
(\ref{117})
into
(\ref{114}),
we know
\begin{align*}
\lim\limits_{n\rightarrow\infty}
\left(
\|v_{n}\|_{W}^{p}-
\mu
\int_{\mathbb{R}^{N}}
\frac{|v_{n}|^{p}}{|x|^{p}}
\mathrm{d}x
\right)
\rightarrow
S_{p,1,\alpha,\mu}
=
\|v\|_{W}^{p}-
\mu
\int_{\mathbb{R}^{N}}
\frac{|v|^{p}}{|x|^{p}}
\mathrm{d}x,
\end{align*}
and
\begin{align*}
\int_{\mathbb{R}^{N}}
\int_{\mathbb{R}^{N}}
\frac{|v(x)|^{p^{*}_{\alpha}}|v(y)|^{p^{*}_{\alpha}}}{|x-y|^{\alpha}}
\mathrm{d}x
\mathrm{d}y
=1.
\end{align*}
Then
$v$
is an extremal.

In addition,
$|v|\in W^{1,p}(\mathbb{R}^{N})$
and
$|\nabla |v||=|\nabla v|$
a.e. in
$\mathbb{R}^{N}$,
therefore,
$|v|$
is also an extremal,
and then there exist non--negative extremals.

Let
$\bar{v}\geqslant0$
be an extremal.
Denote by $\bar{v}_{*}$ the symmetric--decreasing rearrangement
of
$\bar{v}$
(See \cite[Section 3]{Lieb2001}).
From
\cite{Polya1951}
it follows that
\begin{align}\label{19}
\int_{\mathbb{R}^{N}}
|\nabla \bar{v}_{*}|^{p}
\mathrm{d}x
\leqslant
\int_{\mathbb{R}^{N}}
|\nabla \bar{v}|^{p}
\mathrm{d}x.
\end{align}
According to the simplest rearrangement inequality in \cite[Theorem 3.4]{Lieb2001},
we get
\begin{align}\label{20}
\int_{\mathbb{R}^{N}}
\frac{|\bar{v}|^{p}}{|x|^{p}}
\mathrm{d}x
\leqslant
\int_{\mathbb{R}^{N}}
\frac{|\bar{v}_{*}|^{p}}{|x|^{p}}
\mathrm{d}x.
\end{align}
By using Riesz's rearrangement inequality in \cite[Theorem 3.7]{Lieb2001},
we have
\begin{align}\label{21}
\int_{\mathbb{R}^{N}}
\int_{\mathbb{R}^{N}}
\frac{|\bar{v}(x)|^{p^{*}_{\alpha}}|\bar{v}(y)|^{p^{*}_{\alpha}}}{|x-y|^{\alpha}}
\mathrm{d}x
\mathrm{d}y
\leqslant
\int_{\mathbb{R}^{N}}
\int_{\mathbb{R}^{N}}
\frac{|\bar{v}_{*}(x)|^{p^{*}_{\alpha}}|\bar{v}_{*}(y)|^{p^{*}_{\alpha}}}{|x-y|^{\alpha}}
\mathrm{d}x
\mathrm{d}y.
\end{align}
Combining
(\ref{19}),
(\ref{20})
and
(\ref{21}),
and the fact that
$\mu\geqslant0$,
we get that
$\bar{v}_{*}$
is also an extremal,
and then there exist radially symmetric and nonincreasing extremal.
\qed
\section{Proof of Theorem \ref{theorem2}}
For
$p=2$
and
$s\in(0,1)$,
we give
a estimation of extremal function
$u(x)$.
The proof of Theorem \ref{theorem2}
is based on the Coulomb--Sobolev space
$\mathcal{E}^{s,\alpha,2^{*}_{\alpha,s}}(\mathbb{R}^{N})$
and
the endpoint refined Sobolev inequality in Lemma \ref{lemma3}.

\noindent
{\bf Proof of Theorem \ref{theorem2}:}
In this step,
we show some properties of
radially symmetric, nonincreasing and nonnegative function
$u(x)$.
Let
$\bar{\mu}=4^{s}\frac{\Gamma^{2}(\frac{N+2s}{4})}{\Gamma^{2}(\frac{N-2s}{4})}$.
By the definition of extremal
$u$
(see the proof of Theorem \ref{theorem1}),
we know
\begin{align}\label{22}
\int_{\mathbb{R}^{N}}
\int_{\mathbb{R}^{N}}
\frac{|u(x)-u(y)|^{2}}
{|x-y|^{N+2s}}
\mathrm{d}x
\mathrm{d}y
-
\mu
\int_{\mathbb{R}^{N}}
\frac{|u|^{2}}{|x|^{2}}
\mathrm{d}x
=
S_{2,s,\alpha,\mu},
\end{align}
and
\begin{align}\label{23}
\int_{\mathbb{R}^{N}}
\int_{\mathbb{R}^{N}}
\frac{|u(x)|^{2^{*}_{\alpha,s}}|u(y)|^{2^{*}_{\alpha,s}}}{|x-y|^{\alpha}}
\mathrm{d}x
\mathrm{d}y
=1.
\end{align}
Applying
(\ref{22}),
(\ref{23})
and
the definition of
Coulomb--Sobolev space
$\mathcal{E}^{s,\alpha,2^{*}_{\alpha,s}}(\mathbb{R}^{N})$,
we get
$u\in\mathcal{E}^{s,\alpha,2^{*}_{\alpha,s}}(\mathbb{R}^{N})$.

By using
(\ref{22}),
(\ref{23}),
$u\in\mathcal{E}^{s,\alpha,2^{*}_{\alpha,s}}(\mathbb{R}^{N})$
and
Lemma \ref{lemma3},
 we have
\begin{equation}\label{24}
\begin{aligned}
\|u\|_{L^{\frac{2N}{N-2s}}(\mathbb{R}^{N})}
\leqslant&
C_{4}
\left(
\int_{\mathbb{R}^{N}}
\int_{\mathbb{R}^{N}}
\frac{|u(x)-u(y)|^{2}}
{|x-y|^{N+2s}}
\mathrm{d}x
\mathrm{d}y
\right)
^{\frac{(N-\alpha)(N-2s)}{2N(N+2s-\alpha)}}\\
&~~~~~~~~~\cdot
\left(
\int_{\mathbb{R}^{N}}
\int_{\mathbb{R}^{N}}
\frac{|u(x)|^{\frac{2N-\alpha}{N-2s}}|u(y)|^{\frac{2N-\alpha}{N-2s}}}{|x-y|^{\alpha}}
\mathrm{d}x
\mathrm{d}y
\right)
^{\frac{s(N-2s)}{N(N+2s-\alpha)}}\\
=&
C_{4}
\left(
\int_{\mathbb{R}^{N}}
\int_{\mathbb{R}^{N}}
\frac{|u(x)-u(y)|^{2}}
{|x-y|^{N+2s}}
\mathrm{d}x
\mathrm{d}y
\right)
^{\frac{(N-\alpha)(N-2s)}{2N(N+2s-\alpha)}}\\
\leqslant&
C_{4}
\left(
\left(
\frac{\bar{\mu}}{\bar{\mu}-\mu}
\right)
S_{2,s,\alpha,\mu}
\right)
^{\frac{(N-\alpha)(N-2s)}{2N(N+2s-\alpha)}}.
\end{aligned}
\end{equation}
For any
$0<R<\infty$
and
$B(R)=B(0,R)\subset\mathbb{R}^{N}$,
we obtain
\begin{align*}
C_{4}
\left(
\left(
\frac{\bar{\mu}}{\bar{\mu}-\mu}
\right)
S_{2,s,\alpha,\mu}
\right)
^{\frac{(N-\alpha)(N-2s)}{2N(N+2s-\alpha)}}
\geqslant&
\left(
\int_{\mathbb{R}^{N}}
|u(x)|^{\frac{2N}{N-2s}}
\mathrm{d}x
\right)^{\frac{N-2s}{2N}}\\
\geqslant&
\left(
\int_{B(R)}
|u(x)|^{\frac{2N}{N-2s}}
\mathrm{d}x
\right)^{\frac{N-2s}{2N}}\\
\geqslant&
|u(R)|
\omega_{N-1}^{\frac{N-2s}{2N}}
\left(
\int_{0}^{R}
\rho^{N-1}
\mathrm{d}\rho
\right)
^{\frac{N-2s}{2N}}\\
=&
|u(R)|
\left(
\frac{\omega_{N-1}}{N}
\right)
^{\frac{N-2s}{2N}}
R
^{\frac{N-2s}{2}},
\end{align*}
which implies
\begin{align*}
C_{4}
\left(
\left(
\frac{\bar{\mu}}{\bar{\mu}-\mu}
\right)
S_{2,s,\alpha,\mu}
\right)
^{\frac{(N-\alpha)(N-2s)}{2N(N+2s-\alpha)}}
\left(
\frac{N}{\omega_{N-1}}
\right)
^{\frac{N-2s}{2N}}
\frac{1}{
|x|^{\frac{N-2s}{2}}}
\geqslant&
|u(x)|.
\end{align*}
\qed
\section{Proof of Theorem \ref{theorem3}}
For
$p\not=2$
and
$s\in(0,1)$,
we give
a estimation of extremal function
$u(x)$.
From Theorem \ref{theorem1},
we know that
$u(x)$
is a radially symmetric, nonincreasing and nonnegative function.

The proof of Theorem \ref{theorem3}
is different from
Theorem \ref{theorem2}.
The endpoint refined Sobolev inequality in Lemma \ref{lemma3}
is true for
$p=2$.
However,
we don't know that
the endpoint refined Sobolev inequality
is true or not for
$p\not=2$.

\noindent
{\bf Proof of Theorem \ref{theorem3}:}
Let
$\tilde{\mu}=\left(\frac{N-p}{p}\right)^{p}$.
By the definition of extremal
$u$,
we know
\begin{align}\label{25}
\|u\|_{W}^{p}
-
\mu
\int_{\mathbb{R}^{N}}
\frac{|u|^{p}}{|x|^{p}}
\mathrm{d}x
=
S_{p,1,\alpha,\mu},
\end{align}
and
\begin{align}\label{26}
\int_{\mathbb{R}^{N}}
\int_{\mathbb{R}^{N}}
\frac{|u(x)|^{p^{*}_{\alpha}}|u(y)|^{p^{*}_{\alpha}}}{|x-y|^{\alpha}}
\mathrm{d}x
\mathrm{d}y
=1.
\end{align}
For any
$R\in(0,\infty)$
and
$B(R)=B(0,R)\subset\mathbb{R}^{N}$,
by
H\"{o}lder's inequality,
we obtain
\begin{align*}
\left(
\int_{B(R)}
|u|^{p^{*}_{\alpha}}
\mathrm{d}x
\right)^{\frac{1}{p^{*}_{\alpha}}}
\leqslant&
\left[
\left(
\int_{B(R)}
\mathrm{d}x
\right)^{1-\frac{p^{*}_{\alpha}}{p^{*}}}
\left(
\int_{B(R)}
|u|^{p^{*}_{\alpha}\cdot\frac{p^{*}}{p^{*}_{\alpha}}}
\mathrm{d}x
\right)^{\frac{p^{*}_{\alpha}}{p^{*}}}
\right]
^{\frac{1}{p^{*}_{\alpha}}}\\
=&
|B(R)|
^{\frac{1}{p^{*}_{\alpha}}-\frac{1}{p^{*}}}
\left(
\int_{B(R)}
|u|^{p^{*}}
\mathrm{d}x
\right)^{\frac{1}{p^{*}}}\\
\leqslant&
|B(R)|
^{\frac{1}{p^{*}_{\alpha}}-\frac{1}{p^{*}}}
S_{p,1,0,\mu}^{-\frac{1}{p}}\|u\|_{W}\\
\leqslant&
|B(R)|^{\frac{1}{p^{*}_{\alpha}}-\frac{1}{p^{*}}}
S_{p,1,0,\mu}^{-\frac{1}{p}}
\left(
\left(
\frac{\tilde{\mu}}{\tilde{\mu}-\mu}
\right)
S_{p,1,\alpha,\mu}
\right)
^{\frac{1}{p}}
<\infty.
\end{align*}
By Fubini's theorem,
we get
\begin{align*}
(2R)^{-\alpha}
\left(
\int_{B(R)}
|u(x)|^{p^{*}_{\alpha}}
\mathrm{d}x
\right)^{2}
=&
(2R)^{-\alpha}
\int_{B(R)}
|u(x)|^{p^{*}_{\alpha}}
\mathrm{d}x
\int_{B(R)}
|u(y)|^{p^{*}_{\alpha}}
\mathrm{d}y\\
=&
(2R)^{-\alpha}
\int_{B(R)}
\int_{B(R)}
|u(x)|^{p^{*}_{\alpha}}
|u(y)|^{p^{*}_{\alpha}}
\mathrm{d}x
\mathrm{d}y\\
\leqslant&
\int_{B(R)}
\int_{B(R)}
\frac{|u(y)|^{p^{*}_{\alpha}}|u(x)|^{p^{*}_{\alpha}}}{|x-y|^{\alpha}}
\mathrm{d}x
\mathrm{d}y,
\end{align*}
which implies
\begin{align}\label{27}
\left(
\int_{B(R)}
|u(x)|^{p^{*}_{\alpha}}
\mathrm{d}x
\right)^{2}
\leqslant&
(2R)^{\alpha}
\int_{\mathbb{R}^{N}}
\int_{\mathbb{R}^{N}}
\frac{|u(y)|^{p^{*}_{\alpha}}|u(x)|^{p^{*}_{\alpha}}}{|x-y|^{\alpha}}
\mathrm{d}x
\mathrm{d}y.
\end{align}
According to
(\ref{25}),
(\ref{26})
and
(\ref{27}),
we have
\begin{equation*}
\begin{aligned}
1=&
\int_{\mathbb{R}^{N}}
\int_{\mathbb{R}^{N}}
\frac{|u(x)|^{p^{*}_{\alpha}}|u(y)|^{p^{*}_{\alpha}}}{|x-y|^{\alpha}}
\mathrm{d}x
\mathrm{d}y\\
\geqslant&
(2R)^{-\alpha}
\left(
\int_{B(R)}
|u|^{p^{*}_{\alpha}}
\mathrm{d}x
\right)^{2}\\
\geqslant&
(2R)^{-\alpha}
|u(R)|^{2\cdot p^{*}_{\alpha}}
\left(
\omega_{N-1}
\int_{0}^{R}
\rho^{N-1}
\mathrm{d}\rho
\right)^{2}\\
\geqslant&
\frac{\omega^{2}_{N-1}}{2^{\alpha}N^{2}}
|u(R)|^{2\cdot p^{*}_{\alpha}}
R^{2N-\alpha}.
\end{aligned}
\end{equation*}
Then we know
\begin{equation*}
\begin{aligned}
\left(
\frac{2^{\alpha}N^{2}}{\omega^{2}_{N-1}}
\right)^{\frac{1}{2\cdot p^{*}_{\alpha}}}
\frac{1}{R^{\frac{N-p}{p}}}
\geqslant
|u(R)|.
\end{aligned}
\end{equation*}
Hence,
for any
$0<|x|<\infty$,
we obtain
\begin{equation*}
\begin{aligned}
\left(
\frac{2^{\alpha}N^{2}}{\omega^{2}_{N-1}}
\right)^{\frac{1}{2\cdot p^{*}_{\alpha}}}
\frac{1}{|x|^{\frac{N-p}{p}}}
\geqslant
u(x).
\end{aligned}
\end{equation*}
\qed
\section{Conclusions and future works}
The results in this paper set the foundation for the study of a number of questions
related to minimizing problem
$$
S_{p,1,\alpha,\mu}:=
\inf_{u\in W^{1,p}(\mathbb{R}^{N})\setminus\{0\}}
\frac{
\int_{\mathbb{R}^{N}}|\nabla u|^{p}\mathrm{d}x
-
\mu
\int_{\mathbb{R}^{N}}
\frac{|u|^{p}}{|x|^{p}}
\mathrm{d}x}
{\left(
\int_{\mathbb{R}^{N}}
\int_{\mathbb{R}^{N}}
\frac{|u(x)|^{p^{*}_{\alpha}}|u(y)|^{p^{*}_{\alpha}}}{|x-y|^{\alpha}}
\mathrm{d}x
\mathrm{d}y
\right)^{\frac{p}{2\cdot p^{*}_{\alpha}}}},
$$
where
$N\geqslant3$,
$p\in(1,N)$,
$\alpha\in(0,N)$
and
$\mu\in
\left(
0,
\left(
\frac{N-p}{p}
\right)^{p}
\right)$.

During the preparation of the manuscript we faced several problems which are worth to be tackled in forthcoming investigations.
In the sequel,
we shall formulate some of them:

{\bf (a)}
The challenging problems are to prove the rest of Open problem:
the case of
$N\geqslant3$,
$p\in(1,N)$,
$s\in(0,1)$,
$\alpha\in(0,N)$
and
$\mu\in[0,\mathcal{C}_{N,s,p})$,
and
$\mathcal{C}_{N,s,p}$
is defined in
\cite[Theorem 1.1]{Frank2008}.

{\bf (b)}
In \cite{ZhangBL2017},
the authors studied the following minimizing problem:
$$
I_{2,s,\alpha,\mu}(u,v):=
\inf_{
u,v\in W^{s,2}(\mathbb{R}^{N})\setminus\{0\}}
\frac{
\int_{\mathbb{R}^{N}}
\int_{\mathbb{R}^{N}}
\frac{|u(x)-u(y)|^{2}+|v(x)-v(y)|^{2}}{|x-y|^{N+2s}}
\mathrm{d}x
\mathrm{d}y
-
\mu
\int_{\mathbb{R}^{N}}
\left(
\frac{|u|^{2}}{|x|^{2s}}
+
\frac{|v|^{2}}{|x|^{2s}}
\right)
\mathrm{d}x
}
{\left(
\int_{\mathbb{R}^{N}}
\int_{\mathbb{R}^{N}}
\frac{|u(x)|^{2^{*}_{\alpha,s}}|u(y)|^{2^{*}_{\alpha,s}}+|v(x)|^{2^{*}_{\alpha,s}}|v(y)|^{2^{*}_{\alpha,s}}}
{|x-y|^{\alpha}}
\mathrm{d}x
\mathrm{d}y
\right)^{\frac{1}{ 2^{*}_{\alpha}}}},
$$
where
$N\geqslant3$,
$p=2$,
$s\in(0,1)$,
$\mu\in
\left[
0,4^{s}\frac{\Gamma^{2}(\frac{N+2s}{4})}{\Gamma^{2}(\frac{N-2s}{4})}
\right)$
and
$\alpha\in(0,N)$.

It is worth to extend the study of
$I_{2,\textcolor{red}{s},\alpha,\mu}(u,v)$
to the following minimizing problem:
$$
I_{p,s,\alpha,\mu}(u,v):=
\inf_{
u\in W^{s,p}(\mathbb{R}^{N})\setminus\{0\}}
\frac{
\int_{\mathbb{R}^{N}}
\int_{\mathbb{R}^{N}}
\frac{|u(x)-u(y)|^{p}+|v(x)-v(y)|^{p}}{|x-y|^{N+ps}}
\mathrm{d}x
\mathrm{d}y
-
\mu
\int_{\mathbb{R}^{N}}
\left(
\frac{|u|^{p}}{|x|^{ps}}
+
\frac{|v|^{p}}{|x|^{ps}}
\right)
\mathrm{d}x
}
{\left(
\int_{\mathbb{R}^{N}}
\int_{\mathbb{R}^{N}}
\frac{
|u(x)|^{p^{*}_{\alpha,s}}
|u(y)|^{p^{*}_{\alpha,s}}
+
|v(x)|^{p^{*}_{\alpha,s}}
|v(y)|^{p^{*}_{\alpha,s}}
}
{|x-y|^{\alpha}}
\mathrm{d}x
\mathrm{d}y
\right)^{\frac{p}{2\cdot p^{*}_{\alpha,s}}}},
$$
where
$N\geqslant3$,
$p\in(1,N)$,
$s\in(0,1)$,
$\mu\in[0,\mathcal{C}_{N,s,p}$)
and
$\alpha\in(0,N)$.

{\bf (c)}
By using Theorem \ref{theorem1}
and
Lemma
\ref{lemma4},
we could
study
the
Choquard--equation involving two critical nonlinearities:
\begin{equation*}
\begin{aligned}
-\Delta_{p} u
-\mu
\frac{|u|^{p-2}u}{|x|^{p}}
=
\left(
\int_{\mathbb{R}^{N}}
\frac{|u|^{p^{*}_{\alpha}}}{|x-y|^{\alpha}}
\mathrm{d}y
\right)
|u|^{p^{*}_{\alpha}-2}u
+
|u|^{p^{*}-2}u,
\mathrm{~in~}
\mathbb{R}^{N},
\end{aligned}
\end{equation*}
where
$N\geqslant3$,
$p\in(1,N)$,
$\mu\in[0,\left(
\frac{N-p}{p}
\right)^{p}$)
and
$\alpha\in(0,N)$.

{\bf (d)}
Ambrosetti,
Brezis
and
Cerami
\cite{Ambroseiti1994}
proved the existence of infinity many  solutions to the following problem
\begin{equation*}
\begin{aligned}
\begin{cases}
-\Delta u
=
|u|^{2^{*}-2}u
+
\lambda|u|^{q-2}u
&
\mathrm{~in~}
\Omega,\\
 u=0&
\mathrm{~in~}
\partial \Omega,
\end{cases}
\end{aligned}
\end{equation*}
where
$\Omega\subset \mathbb{R}^{N}$
is a smooth bounded domain,
$N\geqslant3$,
$\lambda>0$
and
$q\in(1,2)$.
Garcia
and
Peral
\cite{Garcia1991}
proved the existence of infinity many solutions to following problem
\begin{equation*}
\begin{aligned}
\begin{cases}
-\Delta_{p} u
=
|u|^{p^{*}-2}u
+
\lambda|u|^{q-2}u
&
\mathrm{~in~}
\Omega,\\
 u=0&
\mathrm{~in~}
\partial \Omega,
\end{cases}
\end{aligned}
\end{equation*}
where
$-\Delta_{p}$
is the $p$--Laplacian operator,
$\Omega\subset \mathbb{R}^{N}$
is a smooth bounded domain,
$N\geqslant3$,
$\lambda>0$,
$q\in(1,p)$
and
$p^{*}=\frac{NP}{N-p}$.
Gao and Yang
\cite{Gao2017JMAA}
proved the existence of infinity many solutions to following problem
$$
\begin{cases}
-\Delta u
=
\left(
\int_{\Omega}
\frac{|u|^{2^{*}_{\alpha}}}{|x-y|^{\alpha}}
\mathrm{d}y
\right)
|u|^{2^{*}_{\alpha}-2}u
+
\lambda|u|^{q-2}u
&
\mathrm{~in~}
\Omega,\\
u=0&
\mathrm{~in~}
\mathbb{R}^{N}\setminus \Omega,
\end{cases}
$$
where
$\Omega\subset \mathbb{R}^{N}$
is a bounded domain
with $C^{0,1}$ bounded boundary,
$N\geqslant3$,
$\lambda>0$,
$q\in(1,2)$,
$0<\alpha<N$
and
$2^{*}_{\alpha}=\frac{2N-\alpha}{N-2}$
is the critical Hardy--Littlewood--Sobolev upper exponent.

It is natural to ask:
Does there exist a solution to following problem?
$$
\begin{cases}
-
\Delta_{p} u
=
\left(
\int_{\Omega}
\frac{|u|^{p^{*}_{\alpha}}}{|x-y|^{\alpha}}
\mathrm{d}y
\right)
|u|^{p^{*}_{\alpha}-2}u
+
\lambda|u|^{q-2}u
&
\mathrm{~in~}
\Omega,\\
u=0&
\mathrm{~in~}
\mathbb{R}^{N}\setminus \Omega,
\end{cases}
$$
where
$\Omega\subset \mathbb{R}^{N}$
is a bounded domain
with $C^{0,1}$ bounded boundary,
$N\geqslant3$,
$\lambda>0$,
$p\in(1,N)$,
$q\in(1,p)$
and
$0<\alpha<N$.


\begin{thebibliography}{00}
\bibitem{Ambroseiti1994}
A. Ambrosetti, H. Brezis, G. Cerami,
{\it Combined Effects of Concave and Convex Nonlinearities in Some Elliptic Problems},
J. Funct. Anal. {\bf 122} (1994) 519--543.


\bibitem{O.Alves2016}
C. O. Alves, D. Cassani, C. Tasi, M. Yang,
{\it Existence and concentration of ground state solution for a critical nonlocal Schr\"{o}dinger equation in $\mathbb{R}^{2}$},
J. Differential Equations
{\bf 261} (2016) 1933--1972.

\bibitem{O.Alves2017}
C. O. Alves, F. Gao, M. Squassina, M. Yang,
{\it Singularly Perturbed critical Choquard equations},
J. Differential Equations
{\bf 263} (2017) 3943--3988.

\bibitem{Brezis1983}
H. Brezis, E. Lieb,
{\it A Relation between pointwise convergence of functions and convergence of functional},
Proc. Amer. Math. Soc. 88 (1983) 486--490.

\bibitem{Bellazzini2016}
J. Bellazzini, M. Ghimenti, C. Mercuri, V. Moroz, J. Van Schaftingen,
{\it Sharp Gagliardo--Nirenberg inequalities in fractional	Coulomb--Sobolev spaces},
T. Am. Math. Soc.
(2017);
Doi: https://doi.org/10.1090/tran/7426.


\bibitem{Frank2008}
R. Frank, R. Seiringer,
{\it
Nonlinear ground state representations and sharp Hardy inequalities},
J. Funct. Anal.
{\bf 255}
(2008)
3407--3430.

\bibitem{Garcia1991}
J. Garcia, I. Peral,
{\it Multiplicity of solutions for elliptic problems with critical exponent or with a nonsymmetric term},
T.   Am. Math. Soc. {\bf 323} (1991) 877--895.

\bibitem{Gao2016}
F. Gao, M. Yang,
{\it On the Brezis--Nirenberg type critical problem for nonlinear Choquard equation},
Sci. China Math. (2016);	Doi: 10.1007/s11425--016--9067--5.

\bibitem{Gao2017JMAA}
F. Gao, M. Yang,
{\it On nonlocal Choquard equations with Hardy--Littlewood--Sobolev critical exponetns},
J. Math. Anal. Appl.
{\bf 448} (2017) 1006--1041.

\bibitem{Lieb2001}
E. Lieb, M. Loss,
{\it Analysis, Gradute Studies in Mathematics},
AMS, Providence (2001).

\bibitem{Mukherjee2017Fractional}
T. Mukherjee, K. Sreenadh,
{\it Fractional Choquard equation with critical nonlinearities},
Nonlinear Differ. Equ. Appl.
{\bf 24} (2017) 63.

\bibitem{Moroz2015}
V. Moroz, J. Van Schaftingen,
{\it Groundstates of nonlinear Choquard equations: Hardy--Littlewood--Sobolev critical exponent},
Commun.  Contemp. Math. {\bf 17} 1550005 (2015).

\bibitem{Moroz2016}
V. Moroz, J. Van Schaftingen,
{\it A guide to the Choquard equation},
J. Fix. Point Theory A.
(2016) 1--41.

\bibitem{Nezza2012}
E. Di Nezza, G. Palatucci, E. Valdinoci,
{\it Hitchhiker's guide to the fractional Sobolev spaces},
Bull. Sci. Math. {\bf 229} (2012) 521--573.


\bibitem{Palatucci2014}
G. Palatucci, A. Pisante,
{\it Improved Sobolev embeddings, profile decomposition, and concentration--compactness for fractional Sobolev spaces},
Calc. Var.
{\bf 50} (2014) 799--829.


\bibitem{Pekar1954}
S. Pekar,
{\it Untersuchung\"{u}ber die Elektronentheorie der Kristalle},
Akademie Verlag, Berlin, 1954.

\bibitem{Filippucci2008}
R. Filippucci, P. Pucci, V. R\"{a}dulescu,
{\it Existence and Non--Existence Results for Quasilinear Elliptic Exterior Problems with Nonlinear Boundary Conditions},
Commun.  Part. Diff. Eq. {\bf 33} (2008) 706--717.

\bibitem{Penrose1996}
R. Penrose,
{\it On gravity's role in quantum sstate reduction},
Gen. Relativity Gravitation
{\bf 28} (1996) 581--600.

\bibitem{Polya1951}
G. P\'{o}lya and G. Szeg\"{o},
{\it Isoperimetric Inequalities in Mathematical Physics},
Annals of Mathematics Studies, 27, Princeton University Press, (1951).


\bibitem{Pucci2017}
P. Pucci,
M. Xiang,
B. Zhang,
{\it Existence results for Schrodinger--Choquard--Kirchhoff equations involving the fractional p--Laplacian},
Adv.  Calc.  Var.
(2017) published online.

\bibitem{ZhangBL2017}
L. Wang,
B. Zhang,
H. Zhang,
{\it Fractional Laplacian system involving doubly critical nonlinearities in $\mathbb{R}^{N}$},
Electron. J. Qual. Theo.
{\bf 57} (2017) 1--17.

\bibitem{Willem1996}
M. Willem,
{\it Minimax theorems},
Birkh\"{a}user, Boston, (1996).

\bibitem{Xiang2015}
M. Xiang, Mingqi, B. Zhang, M. Ferrara,
{\it Existence of solutions for Kirchhoff type problem involving the non-local fractional p--Laplacian},
J. Math. Anal. Appl. {\bf 424} (2015) 1021--1041.

\bibitem{Xiang2016}
M. Xiang, Mingqi, B. Zhang, V. R\u{a}dulescu,
{\it Existence of solutions for a bi--nonlocal fractional p--Kirchhoff type problem},
Comput.  Math.  Appl. {\bf 71} (2016) 255--266.

\bibitem{Yang2017}
J. Yang,
F. Wu,
{\it Doubly critical problems involving fractional laplacians in $\mathbb{R}^{N}$},
Adv. Nonlinear Stud.
(2017) published online.
\end{thebibliography}
\end{document}